# Determining full conditional independence by low-order conditioning

DHAFER MALOUCHE

*Unité Signaux et Systèmes, Ecole Nationale d'Ingénieurs de Tunis BP 37, LE BELVEDERE 1002, Tunis, Tunisia. E-mail: dhafer.malouche@essai.rnu.tn*

A concentration graph associated with a random vector is an undirected graph where each vertex corresponds to one random variable in the vector. The absence of an edge between any pair of vertices (or variables) is equivalent to full conditional independence between these two variables given all the other variables. In the multivariate Gaussian case, the absence of an edge corresponds to a zero coefficient in the precision matrix, which is the inverse of the covariance matrix.

It is well known that this concentration graph represents some of the conditional independencies in the distribution of the associated random vector. These conditional independencies correspond to the "separations" or absence of edges in that graph. In this paper we assume that there are no other independencies present in the probability distribution than those represented by the graph. This property is called the perfect Markovianity of the probability distribution with respect to the associated concentration graph. We prove in this paper that this particular concentration graph, the one associated with a perfect Markov distribution, can be determined by only conditioning on a limited number of variables. We demonstrate that this number is equal to the maximum size of the minimal separators in the concentration graph.

## 1. Introduction

Let $(\mathcal{X}, G, \mathcal{F})$ be a triple where $\mathcal{X} = \bigtimes_{\alpha \in V} \mathcal{X}_\alpha$ is a product probability space, $G = (V, E)$ is a graph with a finite set of vertices $V$ and a set of edges $E \subseteq V \times V$ in which a certain separation criteria $\mathcal{C}$ is defined, and $\mathcal{F}$ is a family of probability distribution of random vectors $\mathbf{X} = (X_\alpha, \alpha \in V)'$ with values in $\mathcal{X}$. The triple $(\mathcal{X}, G, \mathcal{F})$ is called a *graphical model* if it satisfies the following property called the *global* Markov property.

Let $A$, $B$ and $S$ be three disjoint subsets where $A$ and $B$ are non-empty. If $S$ separates $A$ and $B$ according to the criteria $\mathcal{C}$ in $G$, denoted by $A \perp_\mathcal{C} B | S$, then the random vectors $\mathbf{X}_A$ and $\mathbf{X}_B$ are independent given $\mathbf{X}_S$, where $\mathbf{X}_A$, $\mathbf{X}_B$ and $\mathbf{X}_S$ are subvectors of $\mathbf{X}$ indexed respectively by the subsets of vertices $A$, $B$ and $S$. So

$$A \perp_{\mathcal{C}} B | S \qquad \text{then } \mathbf{X}_A \perp\!\!\!\perp \mathbf{X}_B | \mathbf{X}_S. \tag{1}$$

Note that the graph $G$ should not contain loops – that is, an edge linking one vertex to itself – and any pair of vertices in $G$ is connected at maximum by one edge, that is, there are no multiple edges between any given pair of vertices.

When the graph $G$ has only undirected edges, that is,

$$(\alpha, \beta) \in E \quad \Longleftrightarrow \quad (\beta, \alpha) \in E,$$

the associated graphical model is called a *concentration graph model* (see Lauritzen (1996)). Dempster (1972) studied concentration graph models for Gaussian distributions under the name of covariance selection or *covariance selection models*. The absence of an edge between a given pair of vertices $(\alpha, \beta)$ in the associated graph indicates that the random variable $X_\alpha$ is independent of $X_\beta$ given all the other variables $\mathbf{X}_{-\alpha\beta} = (X_\gamma, \gamma \neq \alpha, \beta)'$. These models are very well studied, especially the ones corresponding to Gaussian distributions (see Whittaker (1990), Lauritzen (1996), Edwards (2000) and, recently, Letac and Massam (2007) and Rajaratnam *et al.* (2008)). The separation criteria defined on such graphs is a simple separation criteria on undirected graphs: $S \subseteq V$ separates two disjoint non-empty subsets $A$ and $B$ of $V$ if any path joining a vertex in $A$ and another in $B$ intersects $S$.

Other graphical models are represented by graphs with *bi-directed* edges. These models are called *covariance graph models*. The absence of an edge between a given pair of vertices $(\alpha, \beta)$ implies that $X_\alpha$ is *marginally* independent from $X_\beta$, denoted $X_\alpha \perp\!\!\!\perp X_\beta$. The separation criteria in bi-directed graphs can be defined as follows: If $A$, $B$ and $S$ are three disjoint subsets of $V$, where $S$ could be empty, the subset $S$ *separates* $A$ and $B$ in the bi-directed graph $G$ if $V \setminus (A \cup B \cup S)$ separates $A$ and $B$, that is, any path connecting $A$ and $B$ intersects $V \setminus (A \cup B \cup S)$. In this paper this graph will be represented by non-directed edges and will be denoted by $G_0$. So the global Markov property on $G_0$, also called the *covariance global* Markov property, can be defined as follows (see Chaudhuri *et al.* (2007)):

$$\text{If } V \setminus (A \cup B \cup S) \text{ separates } A \text{ and } B \text{ in } G_0, \text{ then } \mathbf{X}_A \perp\!\!\!\perp \mathbf{X}_B | \mathbf{X}_S. \tag{2}$$

Let $P$ be a probability distribution belonging to a certain graphical model $(\mathcal{X}, G, \mathcal{F})$. The probability distribution $P$ is said to be *perfectly Markov* to $G$ if the converse of the *global* Markov property (1) is also satisfied, that is, for any triple of disjoint non-empty subsets $(A, B, S)$ where $S$ is not empty,

$$A \perp_{\mathcal{C}} B | S \quad \Longleftrightarrow \quad \mathbf{X}_A \perp\!\!\!\perp \mathbf{X}_B | \mathbf{X}_S. \tag{3}$$

It was conjectured in Geiger and Pearl (1993) that for any undirected graph $G$ we can find a Gaussian probability distribution $P$ that is *perfectly* Markov to $G$. In the Gaussian case the *perfect* Markovianity assumption is equivalent to the following property: For all non-adjacent vertices $\alpha$ and $\beta$ in $V$ and for all $S \subseteq V \setminus \{\alpha, \beta\}$ and $S \neq \varnothing$,

$$S \text{ separates } \alpha \text{ and } \beta \quad \iff \quad X_\alpha \perp\!\!\!\perp X_\beta | \mathbf{X}_S. \tag{4}$$

In this paper we will consider an undirected graph $G = (V, E)$ and a probability distribution $P$ that is *perfectly Markov* to $G$. Hence, if $\mathbf{X} = (X_\alpha, \alpha \in V)'$ is a random vector with distribution $P$, then $G$ satisfies the following condition

$$(\alpha, \beta) \notin E \quad \iff \quad X_\alpha \perp\!\!\!\perp X_\beta \mid \mathbf{X}_{-\alpha\beta},$$

and if $A$, $B$ and $S$ are three disjoint non-empty subsets of $V$

$$S \text{ separates } A \text{ and } B \text{ in } G \quad \iff \quad \mathbf{X}_A \perp\!\!\!\perp \mathbf{X}_B | \mathbf{X}_S.$$

The aim of this paper is to prove the existence of a relationship between the cardinality of the separators in $G$ and the maximum number of conditioning variables needed to determine the full conditional independencies in $P$. We proceed first by defining a new parameter for undirected graphs, called the "*separability order*". Subsequently we prove that when we condition on a fixed number of variables equal to this *separability order*, the graph that is obtained is exactly the concentration graph.

The paper is organized as follows. Section 2 is devoted to the definition of this parameter and of some properties thereof. In Section 3, we will define a sequence of undirected graphs constructed due to conditional independencies for a given fixed number $k \in \{0, \ldots, |V| - 2\}$. More precisely, we define the $k$-graph $G_k = (V, E_k)$ by

$$(\alpha, \beta) \notin E_k \quad \iff \quad \exists S \subseteq V \setminus \{\alpha, \beta\} \quad \text{such that } |S| = k \text{ and } X_\alpha \perp\!\!\!\perp X_\beta | \mathbf{X}_S. \tag{5}$$

When $k = 0$, the conditional independence given an empty set corresponds to the marginal independence between $X_\alpha$ and $X_\beta$. In the case when $k = 0$, the corresponding $k$-graph is then denoted by $G_0$ and constructed using the pairwise Markov property with respect to bi-directed graphs (see Cox and Wermuth (1996) and Chaudhuri *et al.* (2007)). We mean that

$$\alpha \not\sim_{G_0} \beta \iff X_\alpha \perp\!\!\!\perp X_\beta.$$

The graph $G_0$ is also called a *covariance graph* (see Chaudhuri *et al.* (2007)). Wille and Bühlman (2006) define a graph called a 0–1 graph, which corresponds to a graph with a set of edges equal to $E_0 \cap E_1$. We will show later (see Lemma 7) that this graph is equal to $G_1$. Castello and Roverato (2006) consider so called $k$-partial graphs $G_k^p = (V, E_k^p)$, which are defined as follows:

$$\alpha \not\sim_{G_k} \beta \quad \iff \quad \exists S \subseteq V \setminus \{\alpha, \beta\} \quad \text{such that } |S| \leq k \text{ and } X_\alpha \perp\!\!\!\perp X_\beta | \mathbf{X}_S. \tag{6}$$

Obviously for a fixed $k$, $E_k^p \subseteq E_k$. We will show later (see Lemma 8) that the $k$-partial graph $G_k^p$ is equal to $G_k$. The principle result we prove in this paper (see Theorem 4) is that $E_k \subseteq \cdots \subseteq E_1 \subseteq E_0$ and that $G$ is equal to $G_k$, where $k$ is the *separability order* of $G$. The main assumption of this result is that probability distribution $P$ of the random vector $\mathbf{X}$ is *perfectly* Markov to $G$.

## 2. Separability order

An undirected graph $G = (V, E)$ is a pair of sets where $V$ is the set of vertices and $E$ is the set of edges that is a subset of $V \times V$, where

$$(\alpha, \beta) \in E \iff (\beta, \alpha) \in E.$$

For $\alpha, \beta \in V$, we write $\alpha \sim_G \beta$ when $\alpha$ and $\beta$ are *adjacent* in $G$, that is, $(\alpha, \beta) \in E$. A complete graph is a graph where all the vertices are adjacent, and an empty graph is a graph where the set of edges is empty, that is, $E = \varnothing$. A path between a pair of vertices $(\alpha, \beta)$ is a sequence of distinct vertices $\alpha_0, \alpha_1, \ldots, \alpha_n$ such that $\alpha_0 = \alpha$, $\alpha_n = \beta$ and $\alpha_i \sim_G \alpha_{i+1}$ for all $i = 0, \ldots, (n-1)$.

Let $U \subseteq V$. A subgraph of $G$ induced by $U$ is an undirected graph $G_U = (U, E_U)$ such that $E_U = U \times U \cap E$.

A *connected graph* is a graph $G$ where any pair of vertices can be joined by a path. So a *connected component* of a graph is a subset of $V$ that induces a maximal connected sub-graph of $G$, i.e., $C \subseteq V$ is a connected component of $G$ if $G_C$ is a connected subgraph of $G$ and, for any $\alpha \in V \setminus C$, $\forall \beta \in C$, there is no path between $\alpha$ and $\beta$.

For any non-adjacent vertices $\alpha \not\sim_G \beta$ in a graph $G$ and for any $S \subseteq V \setminus \{\alpha, \beta\}$, we say that $S$ is a *separator* of $\alpha$ and $\beta$ in $G$ if all the paths between $\alpha$ and $\beta$ in $G$ intersect $S$. Consequently, any $S' \supseteq S$ and $S' \subseteq V \setminus \{\alpha, \beta\}$ is also a separator of $\alpha$ and $\beta$. The separator $S$ is called a *minimal separator* of $\alpha$ and $\beta$ if, for any $S' \subseteq S$ and $S' \neq S$, $S'$ cannot be a separator of $\alpha$ and $\beta$. We denote by $\mathrm{ms}_G(\alpha, \beta)$ the set of minimal separators of $\alpha$ and $\beta$ in $G$. It is clear that the set $\mathrm{ms}_G(\alpha, \beta) = \varnothing$ if and only if $\alpha, \beta$ are in two different connected components of $G$.

Let us now give the definition of the *separability order* of an undirected graph:

**Definition 1.** *The separability order of a given graph $G = (V, E)$ is*

$$\mathrm{so}(G) = \max_{\alpha \not\sim_G \beta} \min\{|S|, S \in \mathrm{ms}_G(\alpha, \beta)\} \tag{7}$$

*if $G$ is not complete, and $\mathrm{so}(G) = +\infty$ if $G$ is complete.*

Note that complete graphs have a *separability order* of infinity. Also, empty graphs, that is, graphs with no edges between the vertices of $G$, have a *separability order* equal to zero. Conversely, if $\mathrm{so}(G) = 0$ then either $G$ is composed only of complete connected components or $G$ is an empty graph. We also note that the separability order is purely a graph-theoretic concept.

We now give an example and proceed to prove basic properties of the separability order (see Lemma 1).

**Example 1.** The graph in Figure 1 is an undirected graph containing five vertices. Its *separability order*, $\mathrm{so}(G)$, is equal to 2. We can easily see that

$$\mathrm{ms}_G(1,3) = \{\{2\}\}, \quad \mathrm{ms}_G(2,5) = \{\{3,4\}\},$$
$$\mathrm{ms}_G(1,4) = \{\{2\}\}, \quad \mathrm{ms}_G(1,5) = \{\{3,4\},\{2\}\}.$$

Hence so$(G) = 2$. The degree of the graph $G$, $d(G)$, defined in (10) is equal to 3.

**Lemma 1.** *Let $G = (V, E)$ be an undirected graph with connected components*

$$G_1 = (V_1, E_1), \ldots, G_s = (V_s, E_s),$$

*with $s \geq 2$, and so$(G) = m$. Then*

(i) *$m = 0$ if and only if all the connected components of $G$ are either complete or singletons.*
(ii) *$m = +\infty$ if and only if $G$ is complete.*
(iii) *When $m > 0$, the following properties are satisfied:*
  1. *There exists a pair $(\alpha, \beta)$ of non-adjacent vertices and a minimal separator $S$ of this $\alpha$ and $\beta$ such that $|S| = m$.*
  2. *For any pair of non-adjacent vertices there exists a separator of these two vertices with cardinality equal to $m$.*
(iv) *When $m > 0$, the separability order, so$(G)$, is equal to the maximum separability order among all its non-complete connected components:*

$$\operatorname{so}(G) = \max\{\operatorname{so}(G_l), G_l \ \textit{non-complete}\}. \tag{8}$$

**Proof.** The proof of items (ii) and (iii) follows immediately from Definition 1.

(i) If $m = 0$, this means that for any $\alpha \not\sim_G \beta$ the only separator of these vertices is the empty set. Hence $\alpha$ and $\beta$ belong to different connected components. Moreover, in each connected component of $G$ there are non-adjacent vertices, since $m = 0$. Hence all the connected components of $G$ are either complete or singletons. The converse of this statement follows easily from the definition of the *separability order*.

(iv) Let us define the pairwise separability order of a given pair of vertices $\alpha$ and $\beta$

$$\operatorname{so}(\alpha, \beta | G) = \min\{|S|, S \in \operatorname{ms}_G(\alpha, \beta)\}. \tag{9}$$

Now let $G_1, \ldots, G_l$ be the sequence of non-complete connected components of $G$.

Now so$(\alpha, \beta \mid G) = 0$, if $\alpha \in G_i$ and $\beta \in G_j$ when $i \neq j$ and $i, j \in \{1, \ldots, l\}$.

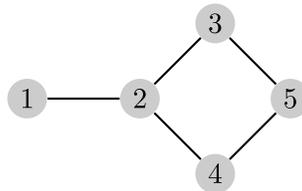

**Figure 1.** so$(G) = 2$ and $d(G) = 3$.

Thus we can focus on the pairwise separability order of pairs within non-complete connected components. Then

$$\begin{aligned}
\mathrm{so}(G) &= \max_{\alpha \not\sim_G \beta} \min\{|S|, S \in \mathrm{ms}_G(\alpha, \beta)\} \\
&= \max_{\alpha, \beta \in V_k, k=1,\ldots,l} \max_{\alpha \not\sim_{G_k} \beta} \min\{|S|, S \in \mathrm{ms}_{G_k}(\alpha, \beta)\} \\
&= \max_{k=1,\ldots,l} \max_{\alpha \not\sim_{G_k} \beta} \min\{|S|, S \in \mathrm{ms}_{G_k}(\alpha, \beta)\} \\
&= \max_{k=1,\ldots,l} \mathrm{so}(G_k).
\end{aligned}$$

$\square$

It is important to note that the *separability order* defined in this paper is exactly equal to the *outer connectivity of the missing edges* defined by Castello and Roverato (2006) for connected graphs. We can also prove that the *separability order* of a non-complete undirected graph $G$ is always smaller than its degree (Lemma 2 below).

**Lemma 2.** *Let $G = (V, E)$ be a non-complete undirected graph; then*

$$\mathrm{so}(G) \leq d(G),$$

*where*

$$d(G) = \max_{\alpha \in V} d(\alpha|G), \tag{10}$$

*where $d(\alpha|G) = |\{\gamma, \alpha \sim_G \gamma\}|$.*

**Proof.** Let $\alpha$ be a vertex in $V$, and let $\mathcal{V}(\alpha|G)$ be the set of vertices adjacent to $\alpha$. So $d(\alpha|G) = |\mathcal{V}(\alpha|G)|$, the degree of $\alpha$ in $G$. Let $\beta$ be a vertex non-adjacent to $\alpha$. It is easy to see that $\mathcal{V}(\alpha|G)$ and $\mathcal{V}(\beta|G)$ are also separators between $\alpha$ and $\beta$. Also it is easy to see that $\mathcal{V}(\alpha|G)$ always contains one minimal separator between $\alpha$ and $\beta$. For example, the set of vertices $\gamma$ denoted by $S(\alpha|G)$ that are simultaneously adjacent to $\alpha$ and belonging to one path between $\alpha$ and $\beta$ is a minimal separator of $\alpha$ and $\beta$. If we suppress one vertex from this $S(\alpha|G)$, this latter set will no longer be a separator between $\alpha$ and $\beta$. The same thing also occurs for $\mathcal{V}(\beta|G)$. Hence

$$\min\{|S|, S \in \mathrm{ms}_G(\alpha, \beta)\} \leq \max(d(\alpha|G), d(\beta|G)).$$

Then,

$$\mathrm{so}(G) = \max_{\alpha \not\sim_G \beta} \{\min\{|S|, S \in \mathrm{ms}_G(\alpha, \beta)\}\} \leq \max_{\alpha \not\sim_G \beta} \max(d(\alpha|G), d(\beta|G))$$

$$= \max_{\alpha \in V} d(\alpha|G) = d(G).$$

So, $\mathrm{so}(G) \leq d(G)$. $\square$

We now define the *degree two* of an undirected graph.

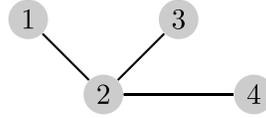

**Figure 2.** $d_2(G) = 1$, $d(G) = 3$, $\text{so}(G) = 1$.

**Definition 2.** *Let $G = (V, E)$ be an undirected graph. The degree two of a vertex in $G$ is defined by,*

$$d_2(\alpha|G) = |\{\gamma, \alpha \sim_G \gamma \text{ and } d(\gamma|G) \geq 2\}|$$

*and the degree two of the graph $G$, $d_2(G)$, is*

$$d_2(G) = \max_{\alpha \in V} d_2(\alpha|G).$$

We give an example to illustrate the *degree two* of a simple undirected graph.

*Example 2.* The graph in Figure 2 has $d_2(G) = 1$, $d(G) = 3$, $\text{so}(G) = 1$.

It is easily seen that in practice the computation of the *separability order* is an NP-complete problem. The degree two of a graph could be a good upper bound for this *separability order*, as this quantity is more easily computable. We prove that $\text{so}(G) \leq d_2(G)$ in Lemma 3 below.

**Lemma 3.** *Let $G = (V, E)$ and $G' = (V, E')$ be two undirected graphs. Then*

  (i) *if $E \subseteq E'$, then $d_2(G) \leq d_2(G')$,*
  (ii) *if $G$ is connected and non-complete $\text{so}(G) \leq d_2(G)$.*

**Proof.**

  (i) First let us define $\mathcal{V}_2(\alpha|G)$ as follows:

  $$\mathcal{V}_2(\alpha|G) = \{\gamma, \alpha \sim_G \gamma \text{ and } d(\gamma|G) \geq 2\}.$$

  Now for any $\alpha \in V$, $\mathcal{V}_2(\alpha|G) \subseteq \mathcal{V}_2(\alpha|G')$. Hence $|\mathcal{V}_2(\alpha|G)| \leq |\mathcal{V}_2(\alpha|G')|$ and thus $d_2(\alpha|G) \leq d_2(\alpha|G')$. This inequality is valid for any $\alpha$ and taking the maximum on $\alpha$ on either side gives $d_2(G) \leq d_2(G')$.
  (ii) If $\text{so}(G) = m$, then, using Lemma 1, part (iii), there exist $\alpha$ and $\beta$ such that $\alpha \not\sim_G \beta$ and a minimal separator $S$ with cardinality $|S| = m$. Now $\mathcal{V}_2(\alpha \mid G)$ contains the set $S(\alpha|G)$, which is a minimal separator between $\alpha$ and $\beta$ (as defined in the proof of Lemma 2). As $S$ is the smallest minimal separator of $\alpha$ and $\beta$,

  $$m \leq |S(\alpha|G)| \leq |\mathcal{V}_2(\alpha|G)| \leq d_2(G)$$

  since $d_2(G) = \max_{\alpha \in V} |\mathcal{V}_2(\alpha|G)|$. Hence $\text{so}(G) \leq d_2(G)$. □

## 3. Concentration graph by low-order conditioning

As before, let $G = (V, E)$ be an undirected graph with the set of vertices $V$ and the set of edges $E$. Let $\mathcal{X} = \times_{\alpha \in V} \mathcal{X}_\alpha$ be a product probability space. The aim of this section is to prove the following result.

**Theorem 4.** *Let $(\mathcal{X}, G, \mathcal{F})$ be a concentration graph model and $P$ a probability distribution belonging to $\mathcal{F}$. Let us consider for any $k \in \{0, \ldots, |V|-2\}$ the undirected graph $G_k = (V, E_k)$ constructed as described in (5):*

$$\alpha \not\sim_{G_k} \beta \iff \exists S \subseteq V \setminus \{\alpha, \beta\} \quad \text{such that } |S| = k \text{ and } X_\alpha \perp\!\!\!\perp X_\beta | \mathbf{X}_S.$$

*Suppose that $P$ is perfectly Markov to $G$ and $\mathrm{so}(G) = m$, then*

$$E = E_m \subseteq E_{m-1} \subseteq \cdots \subseteq E_1.$$

*Furthermore, if $\mathrm{so}(G_0) < |V| - 2$, then $E_1 \subseteq E_0$.*

Theorem 4 will be proved using the following series of lemmas.

**Lemma 5.** *Let $(\mathcal{X}, G, \mathcal{F})$ be a concentration graph model and let $P$ be a probability distribution belonging to $\mathcal{F}$. Suppose that $G$ is connected and non-complete and $\mathrm{so}(G) = m$ where $m > 0$. Suppose also that $P$ is perfectly Markov to $G$. Then $G = G_m$, where $G_m$ is the undirected graph constructed using (5).*

**Proof.** Let $\alpha$ and $\beta$ be two vertices and let us consider a random vector $\mathbf{X} = (X_\alpha, \alpha \in V)'$ with distribution $P$. For any pair $(\alpha, \beta)$ such that $\alpha \not\sim_G \beta$, from Lemma 1(iii), there exists a non-empty subset $S$ with cardinality equal to $m$ that is a separator of $\alpha$ and $\beta$. Using the *global* Markov property with respect to $G$ (see (1)), we can conclude that $X_\alpha \perp\!\!\!\perp X_\beta | \mathbf{X}_S$. Using (5) we conclude that $\alpha \not\sim_{G_m} \beta$. Since this is valid for any pair $(\alpha, \beta)$ we can conclude that $E_m \subseteq E$.

Conversely, suppose that $\alpha \not\sim_{G_m} \beta$; then there exists a separator $S \subseteq V \setminus \{\alpha, \beta\}$ with cardinality $m$ such that $X_\alpha \perp\!\!\!\perp X_\beta | \mathbf{X}_S$. Using the *perfect* Markovianity property we can say that $S$ separates $\alpha$ and $\beta$ in $G$. Thus we can assert that $\alpha \not\sim_G \beta$. Since this argument is valid for any $(\alpha, \beta)$ we can conclude that $E \subseteq E_m$.

We have altogether shown that $E_m \subseteq E$ and $E \subseteq E_m$, hence $E = E_m$, and thus $G = G_m$. □

**Lemma 6.** *Let $(\mathcal{X}, G, \mathcal{F})$ be a concentration graph model and let $P$ be a probability distribution belonging to $\mathcal{F}$. For $m, k \in \mathbb{N}$, let $G_m$ and $G_k$ be two undirected graphs constructed using (5). If $P$ is perfectly Markov to $G$ and $1 \leq m \leq k \leq |V| - 2$ then $E_k \subseteq E_m$.*

**Proof.** Let $\alpha$ and $\beta$ be two vertices. Let us consider a random vector $\mathbf{X} = (X_\alpha, \alpha \in V)'$ with distribution $P$. Suppose that $\alpha \not\sim_{G_m} \beta$, then there exists a separator $S \subseteq V \setminus \{\alpha, \beta\}$ with cardinality $m$ such that $X_\alpha \perp\!\!\!\perp X_\beta | \mathbf{X}_S$. By the *perfect* Markovianity property we

*Determining full conditional independence by low-order conditioning* 1187can conclude that $S$ separates $\alpha$ and $\beta$ in $G$. But since $k \geq m$, we can find a subset $S'$ of $V \setminus \{\alpha, \beta\}$ containing $S$ with cardinality $k$ such that $S'$ is a separator of $\alpha$ and $\beta$ in $G$. Using the *global* Markov property, we determine that $X_\alpha \perp\!\!\!\perp X_\beta | \mathbf{X}_{S'}$. Hence $\alpha \not\sim_{G_k} \beta$.

Since $\alpha \not\sim_{G_m} \beta$ implies that $\alpha \not\sim_{G_k} \beta$ for any pair $(\alpha, \beta)$, we can conclude that $E_k \subseteq E_m$. $\square$

**Lemma 7.** *Let $P$ be a probability distribution in $\mathcal{X}$. Let $G_0 = (V, E_0)$ and $G_1 = (V, E_1)$ be two undirected graphs constructed using (5) for $k = 0, 1$, respectively. If $G_0$ is connected and $\mathrm{so}(G_0) = m_0 < |V| - 2$, then $E_1 \subseteq E_0$.*

**Proof.** Let $\alpha$ and $\beta$ be two vertices and let us consider a random vector $\mathbf{X} = (X_\alpha, \alpha \in V)'$ with distribution $P$. Suppose that $\alpha \not\sim_{G_0} \beta$. By assumption $\mathrm{so}(G_0) = m_0$, hence there exists a subset $S$ of $V \setminus \{\alpha, \beta\}$ for which $|S| = m_0$. Let $\gamma \in V \setminus (S \cup \{\alpha, \beta\})$, which is not empty because $m_0 < |V| - 2$. Then $V \setminus \{\alpha, \beta, \gamma\}$ contains $S$ and so it is a separator of $\alpha$ and $\beta$ in $G_0$. So $\{\gamma\}$ $m$-separates $\{\alpha\}$ and $\{\beta\}$ in $G_0$. Here $G_0$ is seen as an ancestral graph. Using the covariance *global* Markov property (2) with respect to bi-directed graphs (see Cox and Wermuth (1996), Chaudhuri *et al.* (2007), or Drton and Richardson (2008)), we can conclude that $X_\alpha \perp\!\!\!\perp X_\beta | X_\gamma$. Then $\alpha \not\sim_{G_1} \beta$. Hence $E_1 \subseteq E_0$. $\square$

It is easily seen that the results in Lemmas 5–7 lead to the proof of Theorem 4. Castello and Roverato (2006) prove, by also assuming the *perfect* Markovianity, that the concentration graph model $G$ is equal to the $k$-partial graph $G_k^p$ as defined in (6) when $k$ is smaller than the *separability order* of $G$, referred to in that paper as the "outer connectivity of the missing edges". The result in Theorem 4, however, is based on a construction of a sequence of nested graphs. It starts from the covariance graph, that is, $G_0$, the 0-graph, and it becomes stationary and equal to the concentration graph when the number of conditioning variables is equal to the *separability order* of the concentration graph.

In Lemma 8, we show that $k$-partial graphs and $k$-graphs are equal when the *perfect* Markovianity assumption is satisfied. Next, we give a corollary of Theorem 4. In Corollary 9 we give a condition that allows us to determine the last undirected $k$-graph in the sequence of nested graphs obtained due to Theorem 4. This condition is given in term of the *degree two* of the $k$-graphs, not as in Castello and Roverato (2006), where this condition is expressed as a function of the *outer connectivity of connected edges*, a quantity that can be difficult to compute.

**Lemma 8.** *Let $\mathbf{X} = (X_\alpha, \alpha \in V)'$ be a random vector with distribution $P$ belonging to a graphical model generated by an undirected graph $G = (V, E)$. The undirected graphs $G_k^p = (V, E_k^p)$ and $G_k = (V, E_k)$ are respectively the $k$-partial graph and the $k$-graph defined as in (6) and (5). If $P$ is perfectly Markov w.r.t. to $G$ then for any $k \in \{1, \ldots, |V| - 2\}$, we have $E \subseteq E_k = E_k^p$.*

**Proof.** By definition it is easily seen that $\alpha \not\sim_{G_k} \beta$ implies $\alpha \not\sim_{G_k^p} \beta$, hence $E_k \subseteq E_k^p$.

Now let us assume that $\alpha \not\sim_{G_k^p} \beta$, then there exists $S \subseteq V \setminus \{\alpha, \beta\}, |S| \leq k$ such that $X_\alpha \perp\!\!\!\perp X_\beta | \mathbf{X}_S = (X_\gamma, \gamma \in S)'$. If $|S| = k$, the problem is solved. If $|S| < k$, using the *perfect Markovianity* of $P$ we can say that $S$ separates $\alpha$ and $\beta$ in $G$. Then we can construct an $S' \subseteq V \setminus \{\alpha, \beta\}$ with $|S'| = k$, $S' \supseteq S$ such that $S'$ separates $\alpha$ and $\beta$ in $G$. We can now use the *global* Markov property (see (1)), to assert that $X_\alpha \perp\!\!\!\perp X_\beta | \mathbf{X}_{S'}$ and hence $\alpha \not\sim_{G_k} \beta$. We can therefore deduce that $E_k^p \subseteq E_k$. Since $E_k \subseteq E_k^p$ and $E_k^p \subseteq E_k$ we can conclude that $E_k^p = E_k$. The inclusion $E \subseteq E_k$ has already been proved in Theorem 4. □

We can also deduce the following corollary from Theorem 4.

**Corollary 9.** *Let $(\mathcal{X}, G, \mathcal{F})$ be a concentration graphical model such that $G$ is a non-complete connected graph and let $P$ be a probability distribution belonging to $\mathcal{F}$. Let us consider for any $k \in \{0, \ldots, |V| - 2\}$ the undirected graph $G_k = (V, E_k)$ constructed as described in (5). Let us assume that $P$ is perfectly Markov to the graphical model $G$ and $d_2(G) \leq |V| - 2$. Then there exists $k \in \{1, \ldots, |V| - 2\}$ such that*

$$d_2(G_k) \leq k \quad and \quad G = G_k. \tag{11}$$

**Proof.** Let us assume that for all $k \in \{1, \ldots, |V| - 2\}$ that $d_2(G_k) > k$. This implies for example that $d_2(G_{|V|-2}) > |V| - 2$. As the concentration graph $G$ is exactly the $G_{|V|-2}$ we deduce that $d_2(G) > |V| - 2$ which is a contradiction with our assumption, i.e., $d_2(G) \leq |V| - 2$. Hence there exists an integer $k$ such that $d_2(G_k) \leq k$. But $E \subseteq E_k$ and, applying Lemma 3(ii), we deduce that $\text{so}(G) = m \leq d_2(G) \leq d_2(G_k) \leq k$.

Using Theorem 4, as $k \geq m$, we can conclude that $G = G_m = G_k$. □

Corollary 9 can be useful if we wish to determine the concentration graph from a given data set when assuming *perfect* Markovianity. It is sufficient to check the degree two of each estimated $k$-graph.

## 4. Conclusion

In this paper we have proved that a concentration graph model can be determined using a limited number of conditioning variables. The cardinality of this limited subset is determined by looking at the structure of the undirected graph associated with the corresponding distribution *global* Markov property. Certainly the *perfect* Markovianity assumption is also needed for our result to be valid. Our result remains true for both continuous and discrete distributions.

Our result can also be used as a justification of the estimation of graphical models by low-order conditioning such as using the PC algorithm (see Spirtes *et al.* (2000), Kalisch and Bühlmann (2007), Kjærulff and Madsen (2007)), the 0–1 procedure (see Friedman *et al.* (2000) and Wille and Bühlman (2006)), or the *qp*-procedure (see Castello and Roverato (2006)). Practical applications of these procedures are useful when the number of observations are far fewer than the number of variables. We first estimate the sequence of nested graphs $G_k$ starting from $k = 0$ (test on marginal independence between variables) and proceed accordingly. This procedure is terminated when the number of conditioning variables becomes greater than the degree two of the estimated graph $G_k$. In this sense the theory above has tremendous scope for applications.